\documentclass[12pt]{amsart}
\usepackage{amsmath,amssymb,amsthm}
\usepackage{setspace}
\usepackage[colorlinks=true,pagebackref,hyperindex,citecolor=blue,linkcolor=red]{hyperref}
\usepackage{moreverb}
\usepackage{mathrsfs}
\usepackage{url}
\usepackage[all]{xy}
\usepackage{lscape}
\usepackage{tikz}
\usetikzlibrary{arrows.meta,calc}
\usepackage{xcolor}
\usepackage{array}

\newcommand{\ignore}[1]{}

\newtheorem{theorem}{Theorem}[section]

\theoremstyle{definition}
\newtheorem{definition}[theorem]{Definition}
\newtheorem{example}[theorem]{Example}

\theoremstyle{remark}
\newtheorem{remark}[theorem]{Remark}

\numberwithin{equation}{section}

\usepackage{pgfplots}

\pgfplotsset{compat=1.14}

\begin{document}

\title[]{On Maximal order poles of\\ generalized topological zeta functions}
\author[E. Artal-Bartolo]{Enrique Artal Bartolo}
\author[M. Gonz\'alez Villa]{Manuel Gonz\'alez Villa}

\address{Departamento de Matematicas-IUMA, Universidad de Zaragoza, c/ Pedro Cerbuna 12, 50009 Zaragoza, SPAIN}
\address{Centro de Investigaciones Matem\'aticas, A.C., Callej\'on Jalisco s/n, Col. Valenciana\\ C.P. 36023 Guanajuato, Gto., MEXICO} \email{artal@unizar.es, manuel.gonzalez@cimat.mx}

\thanks{First named author is partially supported by
MTM2016-76868-C2-2-P and Grupo ``Álgebra y Geometría'' of
Gobierno de Aragón/Fondo Social Europeo. Second named author is partially supported by
MTM2016-76868-C2-1-P}

\keywords {Topological zeta functions, maximal order poles, plane curves.}

\subjclass[2000]{}

\begin{abstract}
We show some examples of topological  zeta functions associated to an isolated plane curve singular point and an allowed, in the sense of \cite{NV}, differential form that have several poles of order two. This is in contrast to the case of the standard differential form where, as showed by Veys for plane curves \cite{V} and by Nicaise and Xu in general \cite{NX}, there is always at most one pole of order two. 
\end{abstract}

\maketitle

\section{Introduction}

W. Veys determined all poles of the local  topological zeta function of an isolated singular point of a plane curve \cite{V}. In particular, he showed \cite[Theorem 4.2]{V} that there is at most one pole of order two, and, if there is such a pole of order two, it is the largest pole, and its value is the opposite  of the log canonical threshold of the singular point. Later, Veys conjetured that the analogous statements hold for arbitrary dimension \cite[Conjecture 0.2]{LV}, and proved, together with A. Laeremans, the result for polynomials that are non-degenerate with respect to their Newton polyhedron at the origin. They also noticed that a double pole must be of the form $-1/n$ for some $n \in \mathbb{N} \backslash \{0\}$ \cite[Corollary 3.4]{LV}. Finally, the conjecture has recently been established by J. Nicaise and C. Xu \cite[Theorem 3.5]{NX}.

Loeser already proved \cite{L} that if $s_0$ is a pole of order two of topological zeta functions of a singular point of a reduced plane curve, then the monodromy operator of the Milnor fibration of the singular point has  a Jordan block of size two for the eigenvalue exp$(2 \pi i s_0)$. This fact has been generalized, answering a question of C.T.C. Wall \cite{W}, for non reduced plane curves  by A. Melle-Hern\'andez, T. Torelli and W. Veys \cite{MTV}.

Motivated by the remark that the poles of the  topological zeta function determine only a few eigenvalues of the monodromy of the singular point, A. Némethi and W. Veys proposed to study generalized topological zeta functions. Generalized topological zeta functions are associated to a function and an \emph{allowed} differential form. The collection of allowed differential forms for a given funcion $f$ must verify the following three conditions: (1) The standard differential form is allowed. (2) If $s_0$ is a pole of the generalized zeta function associated to the given function and an allowed differential form, then exp$(2 \pi i s_0)$ is an eigenvalue of the monodromy. (3) Any eigenvalue is of the form exp$(2 \pi i s_0)$ for a  pole $s_0$ of  the generalized zeta function associated to some allowed form. 

In this note, we show some examples of generalized  zeta functions associated to an isolated plane curve singular point and an allowed that have several poles of order two. We also discuss on some examples which combinations of order~two poles can appear and find examples where the
double pole is not of the form $-1/n$.

\section{Generalized local topological zeta functions}

Let us define the local topological zeta function associated to an isolated singular  point of a plane curve defined by $f: \mathbb{C}^2 \rightarrow \mathbb{C}$ and a differential form $\omega$.

Let  $\pi : X \rightarrow \mathbb{C}^2$ be an embedded resolution of $f^{-1}(0) \cup {\rm div}(\omega)$. We will consider only holomorphic forms $\omega$ such that an embedded
resolution of $f^{-1}(0)$ is also a resolution of $f^{-1}(0) \cup {\rm div}(\omega)$ and the branching components of both resolutions coincide.
We denote by $E_i$, $i \in S$, the irreducible components (exceptional divisors and strict transforms) 
of the inverse image $\pi^{-1}(f^{-1}(0) \cup {\rm div}(\omega))$. We denote 
by $N_i$ and $\nu_i - 1$ the multiplicities of $E_i$ in the divisor of $\pi^* f$ and $\pi^* \omega$, respectively.  The family $\{(N_i, \nu_i)\}_{i \in S}$ is called the numerical data of the resolution $(X,\pi, f, \omega)$. 
We consider the stratification of $X$ in locally closed subsets given by the subsets  
\[
\mathring{E}_I := 
\left(\bigcap_{i \in I} E_i\right) \setminus\left(\bigcup_{j \not \in I} E_j \right)
\quad \text{for } I \subset S. 
\]
\begin{definition} The (local) topological zeta function of $(f, \omega)$ at $0 \in \mathbb{C}$
is
$$Z_{\text{top}}(f, \omega; s) := \sum_{I \subset S} \chi(\mathring{E}_I \cap \pi^{-1}(0)) \prod_{i \in I} \frac{1}{\nu_i + sN_i} \in \mathbb{Q}(s).$$
\end{definition}

\begin{example}\label{ex:fab}
Let us consider the following two families of cuspidal singular points of complex plane curve.   
\begin{itemize}
\item $g_a=(x-ay^2)^2+y^5$, defines a $(2,5)$ cusp  tangent to $\{x=0\}$ for any $a \in \mathbb{C}$, 
\item $h_b= (y-bx^2)^3+x^7$, defines a $(3,7)$ cusp tangent to $\{y=0\}$ for any $b \in \mathbb{C}$.
\end{itemize}
For any choice $a, b \in \mathbb{C}$ with $a, b \ne 1$, and $a \ne b$, the product
$$
f_{a,b}:= g_{1} \cdot g_{a} \cdot g_{b} \cdot h_1 \cdot h_{-1}
$$
defines an isolated singular point of a plane curve with 5 branches. 
The dual resolution graph and the numerical data of $f_{a,b}$ 
and the differential form
$\omega =  dxdy$ are shown in Figure~\ref{fig:5branches}. Notice that the default numerical data of arrows are $(1,1)$.

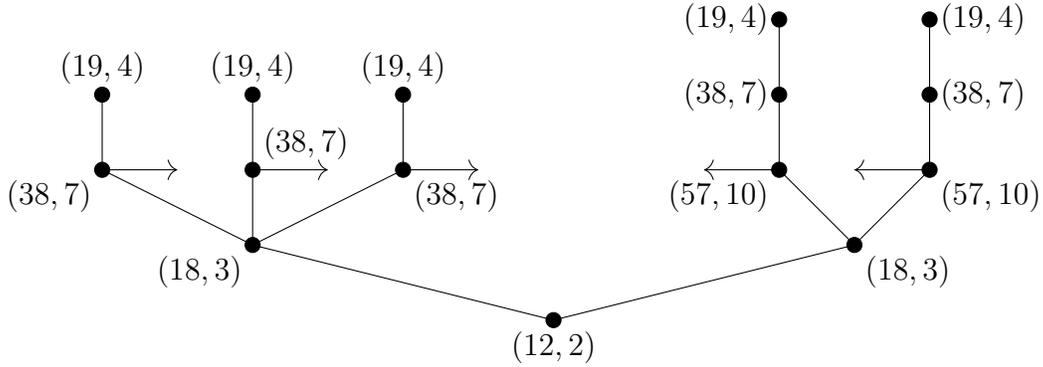
\begin{figure}[ht]
\begin{center}
\begin{tikzpicture}[vertice/.style={draw,circle,fill,minimum size=0.2cm,inner sep=0}]
\coordinate (U) at (0,0);
\coordinate (A1) at (-4,1);
\coordinate (B1) at (4,1);
\coordinate (A11) at (-6,2);
\coordinate (A12) at (-4,2);
\coordinate (A13) at (-2,2);
\coordinate (B11) at (5,2);
\coordinate (B12) at (3,2);
\coordinate (A21) at (-6,3);
\coordinate (A22) at (-4,3);
\coordinate (A23) at (-2,3);
\coordinate (B21) at (5,3);
\coordinate (B22) at (3,3);
\coordinate (B31) at (5,4);
\coordinate (B32) at (3,4);

\node[vertice] at (U) {};
\node[below] at (U) {$(12,2)$};
\node[vertice] at (A1) {};
\node[below left] at (A1) {$(18,3)$};
\node[vertice] at (B1) {};
\node[below right] at (B1) {$(18,3)$};
\node[vertice] at (A11) {};
\node[below left] at (A11) {$(38,7)$};
\node[vertice] at (A12) {};
\node[above right] at (A12) {$(38,7)$};
\node[vertice] at (A13) {};
\node[below right] at (A13) {$(38,7)$};
\node[vertice] at (B11) {};
\node[below right] at (B11) {$(57,10)$};
\node[vertice] at (B12) {};
\node[below left] at (B12) {$(57,10)$};
\node[vertice] at (A21) {};
\node[above] at (A21) {$(19,4)$};
\node[vertice] at (A22) {};
\node[above] at (A22) {$(19,4)$};
\node[vertice] at (A23) {};
\node[above] at (A23) {$(19,4)$};
\node[vertice] at (B21) {};
\node[right] at (B21) {$(38,7)$};
\node[vertice] at (B22) {};
\node[left] at (B22) {$(38,7)$};
\node[vertice] at (B31) {};
\node[right] at (B31) {$(19,4)$};
\node[vertice] at (B32) {};
\node[left] at (B32) {$(19,4)$};

\draw (U)--(A1)--(A11)--(A21);
\draw (A1)--(A12)--(A22);
\draw (A1)--(A13)--(A23);
\draw (U)--(B1)--(B11)--(B21)--(B31);
\draw (B1)--(B12)--(B22)--(B32);
\draw[-{[scale=1.5]>}] (A11)--($(A11)+(1,0)$) ;
\draw[-{[scale=1.5]>}] (A12)--($(A12)+(1,0)$) ;
\draw[-{[scale=1.5]>}] (A13)--($(A13)+(1,0)$) ;
\draw[-{[scale=1.5]>}] (B11)--($(B11)-(1,0)$) ;
\draw[-{[scale=1.5]>}] (B12)--($(B12)+(-1,0)$) ;
\end{tikzpicture}
\caption{Dual graph of the embedded resolution of $f_{a,b}$.}
\label{fig:5branches}
\end{center}
\end{figure}

The Milnor number  of $f_{ab}^{-1}(0)$ at the origin is 188 and the monodromy operator associated to its Milnor fibration has 9 Jordan blocks of size $2$. The eigenvalues (resp. the eigenvalues corresponding to the Jordan blocks of size $2$) are given by the roots  of 
$$
\frac{(t-1)(t^{57}-1)^2(t^{38}-1)^3(t^{18}-1)^3}{(t^{19}-1)^5} \qquad (\text{resp.} \quad \frac{(t^3-1)(t^6-1)(t^2-1)^2}{(t-1)^4}).
$$
The local topological zeta function has a pole ($-1/6$) of order two, due to the subgraph formed by the three vertices with decorations $(12,2)$
and $(18,3)$ and the edges between them,  and it is given by 
$$Z_{\text{top}}(f_{ab},\omega; s)=
\frac{70 +1051 s+5138 s^2+7864 s^3 -1368 s^{4}}{{\left(57 \, s + 10\right)} {\left(38 \, s + 7\right)} {\left(6 \, s + 1\right)}^{2} {\left(s + 1\right)}}
.$$
\end{example}

\section{Examples with several poles of order two}

The following examples show how to produce several poles of order two considering allowed differential forms.

\begin{example}\label{ex:fab1} Let us consider $f_{ab}$ from Example~\ref{ex:fab}
together with the differential form $\omega_1=x^3dxdy$.  Assuming that $a \neq 0$, the dual resolution graph 
is in Figure~\ref{fig:omega1}, showing also the numerical data of $f$ 
and the differential form
$\omega_1$. Again, the default numerical data of (solid) arrows are $(1,1)$; dashed rows correspond
to the strict transform of $(\omega_1)$.
\begin{figure}[ht]
\begin{center}
\begin{tikzpicture}[vertice/.style={draw,circle,fill,minimum size=0.2cm,inner sep=0}]
\coordinate (U) at (0,0);
\coordinate (A1) at (-4,1);
\coordinate (B1) at (4,1);
\coordinate (A11) at (-6,2);
\coordinate (A12) at (-4,2);
\coordinate (A13) at (-2,2);
\coordinate (B11) at (5,2);
\coordinate (B12) at (3,2);
\coordinate (A21) at (-6,3);
\coordinate (A22) at (-4,3);
\coordinate (A23) at (-2,3);
\coordinate (B21) at (5,3);
\coordinate (B22) at (3,3);
\coordinate (B31) at (5,4);
\coordinate (B32) at (3,4);

\node[vertice] at (U) {};
\node[below] at (U) {$(12,5)$};
\node[vertice] at (A1) {};
\node[below left] at (A1) {$(18,9)$};
\node[vertice] at (B1) {};
\node[below right] at (B1) {$(18,6)$};
\node[vertice] at (A11) {};
\node[below left] at (A11) {$(38,19)$};
\node[vertice] at (A12) {};
\node[above right] at (A12) {$(38,19)$};
\node[vertice] at (A13) {};
\node[below right] at (A13) {$(38,19)$};
\node[vertice] at (B11) {};
\node[below right] at (B11) {$(57,19)$};
\node[vertice] at (B12) {};
\node[below left] at (B12) {$(57,19)$};
\node[vertice] at (A21) {};
\node[above] at (A21) {$(19,10)$};
\node[vertice] at (A22) {};
\node[above] at (A22) {$(19,10)$};
\node[vertice] at (A23) {};
\node[above] at (A23) {$(19,10)$};
\node[vertice] at (B21) {};
\node[right] at (B21) {$(38,13)$};
\node[vertice] at (B22) {};
\node[left] at (B22) {$(38,13)$};
\node[vertice] at (B31) {};
\node[right] at (B31) {$(19,7)$};
\node[vertice] at (B32) {};
\node[left] at (B32) {$(19,7)$};

\draw (U)--(A1)--(A11)--(A21);
\draw (A1)--(A12)--(A22);
\draw (A1)--(A13)--(A23);
\draw (U)--(B1)--(B11)--(B21)--(B31);
\draw (B1)--(B12)--(B22)--(B32);
\draw[-{[scale=1.5]>}] (A11)--($(A11)+(1,0)$) ;
\draw[-{[scale=1.5]>}] (A12)--($(A12)+(1,0)$) ;
\draw[-{[scale=1.5]>}] (A13)--($(A13)+(1,0)$) ;
\draw[-{[scale=1.5]>}] (B11)--($(B11)-(1,0)$) ;
\draw[-{[scale=1.5]>}] (B12)--($(B12)+(-1,0)$) ;
\draw[-{[scale=1.5]>},dashed] (A1)--($(A1)+(-2,0)$) node[left] {$(0,4)$} ;
\end{tikzpicture}
\caption{Embedded resolution of $f_{ab}\omega_1$.}
\label{fig:omega1}
\end{center}
\end{figure}
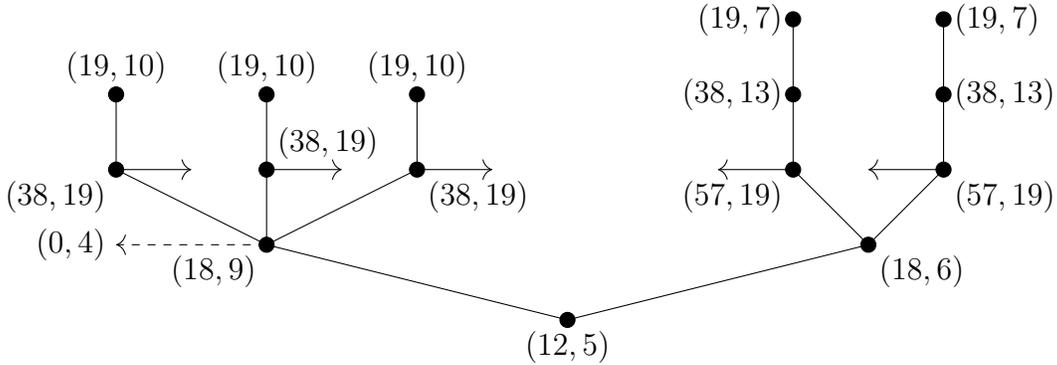

Now, the local topological zeta function has two poles ($-1/2$ and $-1/3$) of order two, and it is given by 
$$
Z_{\text{top}}(f_{ab},\omega_1; s)=
\frac{57+357 \, s+625 \, s^2 + 31\, s^3 - 486 \, s^{4}}{228 \, {\left(3 \, s + 1\right)}^{2} {\left(2 \, s + 1\right)}^{2} {\left(s + 1\right)}}
.$$
The pole $-1/2$ is due to the the vertices with decorations $(18,9)$ and $(38,19)$ and the edges between them. The pole $-1/3$
is due to the the vertices with decorations  $(18,6)$ and $(57,19)$ and the edges between them. 

\end{example}

\begin{example}\label{ex:fab2} Let us consider  $f_{ab}$ from Example~\ref{ex:fab} together with the differential form 
$$\omega_2=x(x-y^2)^2(x-ay^2)^2(x-by^2)^2y^2(y-x^2)^4(y+x^2)^4dxdy.$$
Assuming that $ab \ne 0$, the dual resolution graph and the numerical data of $f_{ab}$ 
and the differential form
$\omega_2$ are in Figure~\ref{fig:omega2}. Again, the default numerical data of (solid) arrows are~$(1,1)$
and dashed arrows correspond to the strict transform of~$(\omega_2)$.

\begin{figure}[ht]
\begin{center}
\begin{tikzpicture}[vertice/.style={draw,circle,fill,minimum size=0.2cm,inner sep=0}]
\coordinate (U) at (0,0);
\coordinate (A1) at (-4,1);
\coordinate (B1) at (4,1);
\coordinate (A11) at (-6,2);
\coordinate (A12) at (-4,2);
\coordinate (A13) at (-2,2);
\coordinate (B11) at (5,2);
\coordinate (B12) at (3,2);
\coordinate (A21) at (-6,3);
\coordinate (A22) at (-4,3);
\coordinate (A23) at (-2,3);
\coordinate (B21) at (5,3);
\coordinate (B22) at (3,3);
\coordinate (B31) at (5,4);
\coordinate (B32) at (3,4);

\node[vertice] at (U) {};
\node[below] at (U) {$(12,19)$};
\node[vertice] at (A1) {};
\node[below left] at (A1) {$(18,27)$};
\node[vertice] at (B1) {};
\node[below right] at (B1) {$(18,30)$};
\node[vertice] at (A11) {};
\node[below left] at (A11) {$(38,57)$};
\node[vertice] at (A12) {};
\node[above right] at (A12) {$(38,57)$};
\node[vertice] at (A13) {};
\node[below right] at (A13) {$(38,57)$};
\node[vertice] at (B11) {};
\node[below right] at (B11) {$(57,95)$};
\node[vertice] at (B12) {};
\node[below left] at (B12) {$(57,95)$};
\node[vertice] at (A21) {};
\node[left] at (A21) {$(19,30)$};
\node[vertice] at (A22) {};
\node[left] at (A22) {$(19,30)$};
\node[vertice] at (A23) {};
\node[right] at (A23) {$(19,30)$};
\node[vertice] at (B21) {};
\node[right] at (B21) {$(38,65)$};
\node[vertice] at (B22) {};
\node[left] at (B22) {$(38,65)$};
\node[vertice] at (B31) {};
\node[right] at (B31) {$(19,35)$};
\node[vertice] at (B32) {};
\node[left] at (B32) {$(19,35)$};

\draw (U)--(A1)--(A11)--(A21);
\draw (A1)--(A12)--(A22);
\draw (A1)--(A13)--(A23);
\draw (U)--(B1)--(B11)--(B21)--(B31);
\draw (B1)--(B12)--(B22)--(B32);
\draw[-{[scale=1.5]>}] (A11)--($(A11)+(1,0)$) ;
\draw[-{[scale=1.5]>}] (A12)--($(A12)+(1,0)$) ;
\draw[-{[scale=1.5]>}] (A13)--($(A13)+(1,0)$) ;
\draw[-{[scale=1.5]>}] (B11)--($(B11)-(1,0)$) ;
\draw[-{[scale=1.5]>}] (B12)--($(B12)+(-1,0)$) ;
\draw[-{[scale=1.5]>},dashed] (A1)--($(A1)+(-2,0)$) node[left] {$(0,2)$} ;
\draw[-{[scale=1.5]>},dashed] (B1)--($(B1)+(2,0)$) node[right] {$(0,3)$} ;
\draw[-{[scale=1.5]>},dashed] (A21)--($(A21)+(0,1)$) node[left] {$(0,3)$} ;
\draw[-{[scale=1.5]>},dashed] (A22)--($(A22)+(0,1)$) node[left] {$(0,3)$} ;
\draw[-{[scale=1.5]>},dashed] (A23)--($(A23)+(0,1)$) node[right] {$(0,3)$} ;
\draw[-{[scale=1.5]>},dashed] (B31)--($(B31)+(0,1)$) node[left] {$(0,5)$} ;
\draw[-{[scale=1.5]>},dashed] (B32)--($(B32)+(0,1)$) node[left] {$(0,5)$} ;
\end{tikzpicture}
\caption{Embedded resolution of $f_{ab}\omega_2$.}
\label{fig:omega2}
\end{center}
\end{figure}
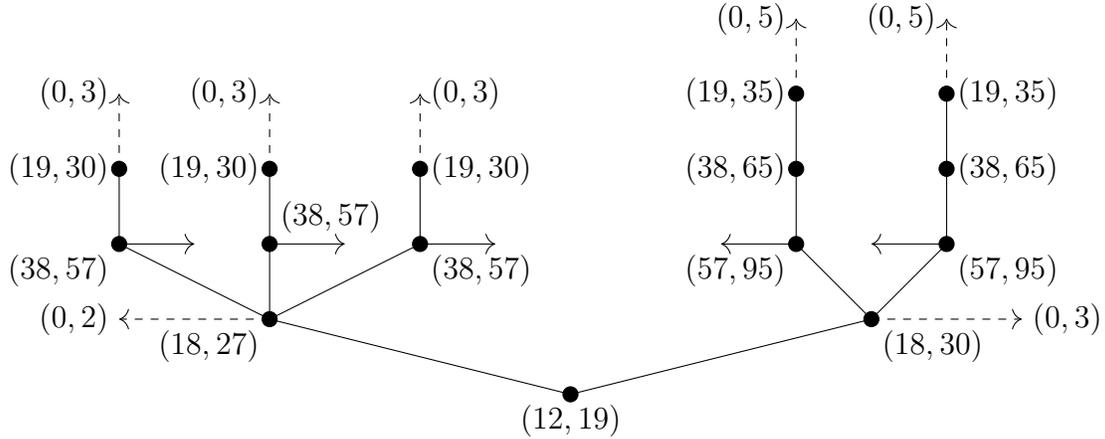

In this example, the local topological zeta function has two poles ($-3/2$ and $-5/3$) of order two,  and it is given by 
$$Z_{\text{top}}(f_{ab},\omega_2; s)=
-\frac{16734 \, s^{4} + 88541 \, s^{3} + 168881 \, s^{2} + 134709 \, s + 36195}{1710 \, {\left(3 \, s + 5\right)}^{2} {\left(2 \, s + 3\right)}^{2} {\left(s + 1\right)}}
.$$
The pole $-3/2$ is due to the the vertices with decorations $(18,27)$ and $(38,57)$ and the edges between them. The pole $-5/3$
is due to the the vertices with decorations  $(18,30)$ and $(57,95)$ and the edges between them. 

\end{example}

\begin{remark}It is not hard to prove that it is impossible
to find any pair of double poles. First,  if the pole attached to the first exceptional component
is double, no other double pole exists. This is the case of the poles of the form $-\frac{1}{6} - k$, and $-\frac{5}{6}-k$, for some $k \in \mathbb{Z}_{>0}$. Notice that $-\frac{1}{6}$ (resp. $-\frac{5}{6}$) is attained with the standard differential form (resp.  with the form $\omega_3=(xy)^4dxdy$), see Examples \ref{ex:fab}, and \ref{ex:fab3}.  Moreover, it is not possible to simultaneously achieve double poles of the forms $-\frac{1}{3}-k$, and  $-\frac{2}{3}-k$, because $3$ does not divide $\gcd(18,38)$. The other combinations are attained in Examples \ref{ex:fab1}, and \ref{ex:fab2}. However, not all choices for $k$ are possible. For instance, one can check that double poles $-\frac{2}{3}$ and $-\frac{3}{2}$
cannot happen simultaneously.
\end{remark}                            

\section{Poles distinct from \texorpdfstring{$-1/n$}{-1/n}}

The following examples show how to produce  poles of order two distinct from $-1/n$ considering allowed differential forms.

\begin{example}\label{ex:fab3} Let us consider consider $f_{ab}$ from Example~\ref{ex:fab} together with the differential form $\omega_3=(xy)^4dxdy$. The dual resolution graph and the numerical data of $f$ 
and the differential form
$\omega_2$ are in Figure~\ref{fig:omega3}. Again, the default numerical data of (solid) arrows are~$(1,1)$
and dashed arrows correspond to the strict transform of~$(\omega_3)$.

\begin{figure}[ht]
\begin{center}
\begin{tikzpicture}[vertice/.style={draw,circle,fill,minimum size=0.2cm,inner sep=0}]
\coordinate (U) at (0,0);
\coordinate (A1) at (-4,1);
\coordinate (B1) at (4,1);
\coordinate (A11) at (-6,2);
\coordinate (A12) at (-4,2);
\coordinate (A13) at (-2,2);
\coordinate (B11) at (5,2);
\coordinate (B12) at (3,2);
\coordinate (A21) at (-6,3);
\coordinate (A22) at (-4,3);
\coordinate (A23) at (-2,3);
\coordinate (B21) at (5,3);
\coordinate (B22) at (3,3);
\coordinate (B31) at (5,4);
\coordinate (B32) at (3,4);

\node[vertice] at (U) {};
\node[below] at (U) {$(12,10)$};
\node[vertice] at (A1) {};
\node[below left] at (A1) {$(18,15)$};
\node[vertice] at (B1) {};
\node[below right] at (B1) {$(18,15)$};
\node[vertice] at (A11) {};
\node[below left] at (A11) {$(38,31)$};
\node[vertice] at (A12) {};
\node[above right] at (A12) {$(38,31)$};
\node[vertice] at (A13) {};
\node[below right] at (A13) {$(38,31)$};
\node[vertice] at (B11) {};
\node[below right] at (B11) {$(57,46)$};
\node[vertice] at (B12) {};
\node[below left] at (B12) {$(57,46)$};
\node[vertice] at (A21) {};
\node[above] at (A21) {$(19,16)$};
\node[vertice] at (A22) {};
\node[above] at (A22) {$(19,16)$};
\node[vertice] at (A23) {};
\node[above] at (A23) {$(19,16)$};
\node[vertice] at (B21) {};
\node[right] at (B21) {$(38,31)$};
\node[vertice] at (B22) {};
\node[left] at (B22) {$(38,31)$};
\node[vertice] at (B31) {};
\node[right] at (B31) {$(19,16)$};
\node[vertice] at (B32) {};
\node[left] at (B32) {$(19,16)$};

\draw (U)--(A1)--(A11)--(A21);
\draw (A1)--(A12)--(A22);
\draw (A1)--(A13)--(A23);
\draw (U)--(B1)--(B11)--(B21)--(B31);
\draw (B1)--(B12)--(B22)--(B32);
\draw[-{[scale=1.5]>}] (A11)--($(A11)+(1,0)$) ;
\draw[-{[scale=1.5]>}] (A12)--($(A12)+(1,0)$) ;
\draw[-{[scale=1.5]>}] (A13)--($(A13)+(1,0)$) ;
\draw[-{[scale=1.5]>}] (B11)--($(B11)-(1,0)$) ;
\draw[-{[scale=1.5]>}] (B12)--($(B12)+(-1,0)$) ;
\draw[-{[scale=1.5]>},dashed] (A1)--($(A1)+(-2,0)$) node[left] {$(0,5)$} ;
\draw[-{[scale=1.5]>},dashed] (B1)--($(B1)+(2,0)$) node[right] {$(0,5)$} ;
\end{tikzpicture}
\caption{Embedded resolution of $f_{ab}\omega_3$.}
\label{fig:omega3}
\end{center}
\end{figure}
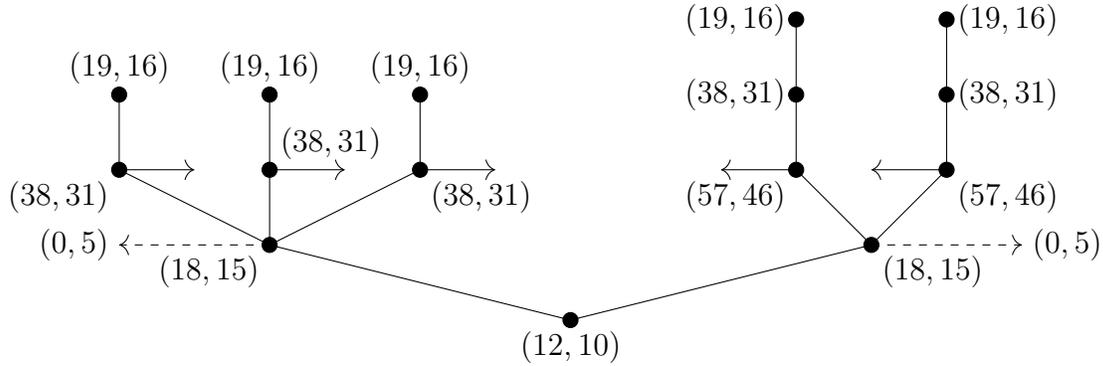

The local topological zeta function has a pole ($-5/6$) of order two,  and it is given by 
$$Z_{\text{top}}(f_{ab},\omega_3; s)=
\frac{-41496 \, s^4 - 91616 \, s^3 - 55926 \, s^2 + 1559 \, s + 7130}{5 \, {\left( 6 \, s + 5 \right)}^2 {\left(38 \, s + 31 \right)} {\left(57 \, s + 46 \right)} {\left( s + 1 \right)}}.$$

\end{example}

\begin{example} Let us consider $g_{p,q}=(y^p+x^q)(y^q+x^p)$, $1<p<q$, $\gcd(p,q)=1$, $\omega=(xy)^{a-1}dxdy$, $a\geq 1$.
An embedded resolution of $(g_{p,q},\omega)$ has too many irreducible components, but in
\cite[Example 4.6]{AMO} we have a $\mathbf{Q}$-resolution with few strata,
see~Figure~\ref{fig:acampo-gen}, where $[\bullet]$ stands for the order of the group associated to the $0$-dimensional stratum, see the explanation below.

\begin{figure}[ht]
\begin{center}
\begin{tikzpicture}[scale=1.5,vertice/.style={draw,circle,fill,minimum size=0.15cm,inner sep=0}]

\draw ($1.25*(0,0)-.25*(-1,1)$) -- ($-.25*(0,0)+1.25*(-1,1)$);

\draw ($1.25*(0,0)-.25*(1,1)$) -- ($-.25*(0,0)+1.25*(1,1)$);

\draw[->] ($(-.5,.5)-.25*(1,1)$) -- ($(-.5,.5)+.25*(1,1)$) node[pos=-.4] {$\boxed{1+s}$};
\draw[->] ($(.5,.5)-.25*(-1,1)$) -- ($(.5,.5)+.25*(-1,1)$) node[above right,pos=.0] {$\boxed{1+s}$};
\draw[->,dashed] ($(-1,1)-.25*(1,1)$) -- ($(-1,1)+.25*(1,1)$) node[pos=-.4] {$\boxed{a}$};
\draw[->,dashed] ($(1,1)-.25*(-1,1)$) -- ($(1,1)+.25*(-1,1)$) node[above,pos=1] {$\boxed{a}$};

\node[vertice] at (0,0) {};
\node[right] at (.1,0) {$[q^2-p^2]$};
\node[vertice] at (1,1) {};
\node[left=2pt] at (1,1) {$[p]$};
\node[vertice] at (-1,1) {};
\node[right=2pt] at (-1,1) {$[p]$};
\node[left] at  ($-.25*(0,0)+1.25*(-1,1)$) {$\boxed{(p+q)(a+ps)}$};
\node[right] at  ($-.25*(0,0)+1.25*(1,1)$) {$\boxed{(p+q)(a+ps)}$};
\end{tikzpicture}
\caption{$\mathbf{Q}$-resolution of $(y^p+x^q)(y^q+x^p)$.}
\label{fig:acampo-gen}
\end{center}
\end{figure}
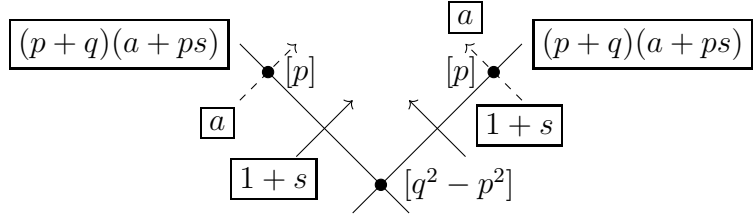

The concept of an embedded $\mathbf{Q}$-resolution
$\pi : X \rightarrow \mathbb{C}^2$ was introduced in~\cite{Qres}. For dimension~$2$, we mean that a finite
set of points in the exceptional divisor may be quotient singularities for some cyclic group; the preimage
of the total divisor satisfies a natural condition of $\mathbb{Q}$-normal crossings. 
In this case there are also some strata $p_1^i,\dots,p_i^{n_i}$ 
such that $(X,p_j^i)\cong(\mathbb{C}^2,0)/G_j^i$ where 
$G_j^i$ is a group of order $m_i^j>1$ with a small action on~$\mathbb{C}^2$. The 
one-dimensional strata will be of the form $\mathring{E}_i=E_i\setminus\left(\bigcup_{j\neq i} E_j\cup\{p_1^1,\dots,p_{n_i}^i\}\right)$. The intersection of two components $E_i\cap E_j$ has also associated an order $m_{i,j}\geq 1$.
Following~\cite{V1} we have:
\[
Z_{\text{top}}(f, \omega; s) = \sum_{i\in I}  \frac{\chi(\mathring{E}_i)+\sum_{j=1}^{n_i}m_i^j}{\nu_i + sN_i}+
\sum_{E_i\cap E_j\neq\emptyset} \frac{m_{i,j}}{(\nu_i + sN_i)(\nu_j + sN_j)} \in \mathbb{Q}(s).
\]

Applying the formula we obtain
\[
Z_{\text{top}}(g_{p,q}, \omega; s)=
\frac{q-p}{q+p}\frac{1}{(a+p s)^2}+
\frac{2}{{\left(a - p\right)} {\left(q+p\right)}}
\left(\frac{1}{1+s}+\frac{{\left(a p - p^{2} - 1\right)} a}{a+p s}\right).
\]
In particular, for $a\not\equiv0\mod{p}$, we obtain that $-\frac{a}{p}$
is a double pole.\end{example}

                                      
\providecommand{\bysame}{\leavevmode\hbox to3em{\hrulefill}\thinspace}
\providecommand{\MR}{\relax\ifhmode\unskip\space\fi MR }
\providecommand{\MRhref}[2]{%
  \href{http://www.ams.org/mathscinet-getitem?mr=#1}{#2}
}
\providecommand{\href}[2]{#2}

\end{document}